# Digital Technologies In The Early Primary School Classroom

Nathalie Sinclair
Simone Fraser University

Anna Baccaglini-Frank
University of Modena and Reggio Emilia (Italy)

**Introduction**

Papert's (1980) work with Turtle Geometry offered an early and provocative vision of how digital technologies could be used with young learners. Since then, research on digital technology use has focused on the middle and high school levels (notable exceptions include Sarama & Clements, 2002; Hoyles, Noss & Adamson, 2002). Given the increasing diversity of digital technologies, and their varied underlying pedagogical goals and design choices, Clements' (2002) claim that "there is no single effect of the computer on mathematics achievement" (p. 174) is as true now as it was a decade ago. However, many advances have been made in better articulating the range of design choices that are possible, their potential effect on the cognitive and affective dimensions on mathematics learning, and their varying demands on the teacher. The aim of this chapter is to summarise the research literature on the use of digital technologies in the teaching and learning of mathematics at the k-2 level. In particular, we focus on literature that contributes to our understanding of *how* the use of digital technologies affects and changes the teaching and learning of mathematics—that is, how different affordances and design choices impact on the way teachers and learners interact and express themselves mathematically. By digital technologies we refer to a range of tools including multi-purpose computer-based software programs, web-based applets, virtual manipulatives, programming languages, CD-ROMs, games, calculators, touchscreen applications and interactive whiteboards. The distinction between these various types is not always evident [1] and, indeed, one goal of this chapter is to provide useful distinguishing features of these various technologies in order to help educators better evaluate and choose amongst them.

We begin by outlining some of the major theoretical developments that are shaping the way researchers are studying the use of digital technologies; we hope that some of these developments, which originate in research conducted for the middle and high school grades, can inform research at the younger grade levels, thus building on decades-old insights and constructs. We then present an overview of research related first to two content areas of the primary school curriculum—number sense and geometry—and second to a mix of content areas all approached through the use of programming languages. Where possible, we try to describe the particular affordances of the digital technology involved, that is, the kinds of interactions that can be performed, acknowledging that intended affordances may not always be perceived as possible by users [2]. We are also aware of the fact that many of the tools we describe may quickly disappear, to be replaced by new interpretations or available on new platforms. We have thus tried to focus attention on the design principles that may have relevance beyond specific examples. At the end of the chapter we discuss several themes that emerge from our survey of the literature and recommend future research directions.

**Historical and theoretical perspectives**

The use of digital technologies in the early grades has traditionally encountered opposition by those concerned that children at this age need tactile, concrete experiences. Indeed, the k-2 classroom has long featured the use of physical manipulatives, with both researchers and teachers acknowledging their importance (Sowell, 1989). This presence of a rich set of resources in the classroom may in fact make it easier for digital technologies to be integrated, in comparison to the higher grades where the technologies of paper-and-pencil usually prevail. Indeed, over the past decade, several researchers have argued for the appropriateness and benefit of using "virtual manipulatives" (VMs) in the early grades, which build on the familiarity of physical ones, but which may also provide a range of added affordances (Bolyard et al., 2010, Moyer-Packenham, 2010; Moyer-Packenham & Suh, 2012; Moyer-Packenham et al., 2013). These researchers have questioned the assumption that "concrete" tools are more appropriate for young children and have argued that physical manipulatives are limited in their ability to promote both mathematical actions and reflections on these actions (Sarama & Clements, 2009). Sarama and Clements authors point specifically to a VMs potential for supporting the development of *integrated-concrete* knowledge, which interconnects knowledge of physical objects, actions on these objects and symbolic representations of these objects and actions. They offer seven hypothesized, interrelated affordances that have been ratified by an admittedly small amount of existing research: bringing mathematical ideas and action to conscious awareness; encouraging and facilitating complete, precise explanations; supporting mental "actions on objects"; changing the very nature of the manipulative; symbolising mathematical concepts; linking the concrete and the symbolic with feedback; and, recording and replaying student actions. Thus, one way to approach the design and evaluation of particular tools is to see whether these affordances are present in the tool in a way that is relevant to the mathematical concept under investigation and accessible to both teachers and learners. We will use these affordances as a way of describing and contrasting the various VMs presented in the next sections.

These affordances, of course, are not unique to VMs and offer a compelling set of macro-level goals for digital technology design. However, they say little about the forms of interaction that different digital technologies might offer. Sedig and Sumner (2006) propose a framework that distinguishes three forms interactions that are "based on three fundamental, root metaphors derived from the way in which humans use their bodies to interact with the external world" (p. 9): conversing, manipulating and navigating. Conversing interactions are ones in which learners issue an input action, which can be done through, for example, procedure-based programming languages, text-based menus and pen-based gestures. Conversing interactions are usually *discrete* in that there is a separation between the learner's actions and the computer-based reaction. For example, clicking a button on the screen that rotates a shape is an example of a discrete conversing interaction. Manipulating interactions involve touching, handling or grasping element(s) on a screen through selecting, dragging and moving. Such interactions are usually considered more tangible than conversing ones in that learners can "reach their hand" into the screen to handle the objects. They are often *continuous* in contrast to *discrete*, as exemplified by the interaction of a learner dragging the vertex of a shape to rotate it, so that cause and effect are observed simultaneously [3]. Navigating interactions involve moving on, over, or through the screen. The majority of digital technologies researched at the k-2 level focus primarily on manipulating and conversing interactions. This may be in part due to the number-focused nature of most of these technologies since navigating interactions are more associated with spatial reasoning. All three interactions can have either a *direct* or *indirect* "focus", this distinction being based on whether the learner is directly interacting with a screen object or interacting with it through an intermediary representation. When a learner is rotating a triangle by dragging one of its vertices, the interaction is direct, but if the learner is dragging a dial in order to rotate the triangle, the interaction is indirect.

Goodwin and Highfield (2013) offer a somewhat different characterization of digital technologies, which focuses more on their constraints and underlying pedagogy: instructive, manipulable and constructive. Instructive digital technologies tend to promote procedural learning, relying on evaluative feedback and repetitive interactions with imposed representations. Manipulable digital technologies enabled the imposed representations to be manipulated so as to engage students in discovery and experimentation. There is much overlap between this category and Sedig and Sumner's manipulating form of interaction, though the former carries with it particular pedagogical goals that the latter does not assume. Finally, constructive digital technologies are ones in which learners create their own representations, which are often the goal of the activity, thereby promoting mathematical modeling and what Noss and Hoyles (1996) characterize as expressive uses of technology. Goodwin and Highfield argue that while instructive technologies may be well-suited for procedural learning, manipulable and constructive technologies better support conceptual learning.

While some digital technologies fit neatly into one particular category (of each of these tripartite characterizations), many will belong to more than one category. However, each characterisation provides a way of comparing the constraints and affordances of different digital technologies, which may guide the choice of a specific digital technology for a particular topic and/or grade level. However, educators must also make larger-scale decisions that involve choosing appropriate digital technologies for a wide range of topics across several grade levels. Is it preferable to promote one category of digital technology or to have a mix of forms of interactions (conversing, manipulating and navigating) and of constraints (*instructive*, *manipulable* and *constructive*)? Although the research has little to say about such a question, Goldenberg (2000) has argued for the "fluency principle", which states that "[l]earning a few good tools well enough to use them knowledgeably, intelligently, mathematically, confidently, and appropriately in solving otherwise difficult problems makes a genuine contribution to a student's mathematical education" (p. 7). A "good" tool might offer a variety of forms of interactions and even enable different kinds of constraints, while also having a consistency that more easily enables teachers and learners to perceive important affordances. Unfortunately, no longitudinal research exploring the effect of long-term use of particular "good" technologies currently exists.

While Sedig and Sumner's framework says very little about *how* mathematical learning happens, several other theories in mathematics education have been proposed with that purpose in mind. These include the instrumental genesis approach, which is primarily concerned with the process of how a computer tool becomes for learners an instrument to learn and do mathematics with (Artigue, 2002); it attends to the way affordances are perceived both through increased experience with the tool and through the problems that the tool enables solving. This theory is specifically devoted to studying technology-based interactions and does not get used by researchers working outside this domain. Very few studies at the primary school level draw on this theory, perhaps because of the nature of computer tools at this age level, which do not require a significant instrumentation process because of their ease-of-use. However, the expanded notions of instrumental orchestration (Trouche, 2004; Drijvers, 2012) and documentational genesis (Gueudet & Trouche, 2008), both of which focus on the work of the teacher in a computer-based classroom, have been used to show, for example, the specific strategies (knowing as orchestration types) that kindergarten teachers use to manage heterogeneity and lack of reading ability at this school level (Gueudet, Bueno & Poisard, 2013).

Another important theoretical approach is that of semiotic mediation (Bartolini-Bussi & Mariotti, 2008), which has its roots in Vygotsky's work approach and attends to the specific ways in which tools (including digital ones) are transformed by learners into mathematical concepts through a process of internalisation. Although focused on tool-use, and developed by researchers working with digital technologies, this theory also concerns mathematical learning more generally. This approach enables researchers to focus on the specific actions that certain tools enable (such as dragging and tracing), and on the types of signs they can give rise to. Similar in its Vygotskian origins, activity theory has also been used in the context of research on technology-based teaching and learning. Ladel and Kortenkamp (2012) propose a specific version of it, which they call artifact-centric activity theory (ACAT), and which they developed specifically for the primary school level, with touchscreen technologies in mind. This approach, more than the others, emphasizes how tools radically change the way learners act and think, thus moving away from a view that tools are discardable crutches that merely scaffold the learning of mathematics. Ladel and Kortenkamp also use it in the very design process of their touchscreen digital technologies.

Another Vygotskian-inspired framework that has been used to investigate technology-based student learning is Sfard's communicational approach, which takes thinking to be communicating and thus learning to be a change in one's discourse (Sfard, 2008). Changes in discourse can involve different uses of particular words and gestures; they can also be based on or produce different visual mediators and different "routines" for identifying shapes or describing quantity. Due to their highly visual and often temporal nature, digital technologies quite frequently offer unique visual mediators, thus inviting different ways of describing and comparing mathematical objects and relationships. They also give rise to new ways of thinking that may conflict with the established discourse of formal mathematics (which tends to be static and alphanumeric); Sfard's approach can help draw attention to the potential communication conflicts that may thus arise, especially as teachers attempt to transition between digital and text-based resources (see Sinclair & Yurita, 2008).

Much of the research on the use of digital technologies has also been informed by theories of embodied cognition. Papert's notion of "body syntonicity" can be seen as an early precursor to the now widely-shared recognition of the important role that the body plays in mathematical meaning-making. While there are a range of assumptions about the relationship between the body and the mind—with dualist perspectives seeing the body as an important and sometimes necessary scaffolding for the development of mathematical schemas and concept and the monist perspective seeing the body itself as doing the knowing—there is growing consensus that particular kinds of bodily engagement can support mathematics learning. Research focused on the particular ways in which digital technologies can enable and promote bodily engagement highlights the precise and temporal actions that these technologies afford, which enable learners to move in mathematically-relevant ways (Nemirovsky, Kelton & Rhodehamel, 2013; Robutti, 2006; Sinclair, de Freitas & Ferrara, 2013). Indeed, the three basic metaphors mentioned earlier—conversing, manipulating and navigating—provide ways for the speaking, hearing, touching and seeing body to interact with mathematical objects. This multimodal kind of interaction seems particularly important at the primary school level, where children's communication is much less language-based—we return to this point in our discussion of new, touchscreen technologies.

More recently, research focused on learning with media (such as television or videos) suggests that joint engagement—which involves parental mediation—can provide powerful additional affordances for learning beyond what is found with technology use alone (Moorthy et al., 2013; Stevens & Penuel, 2010). This work has been extended to the context of digital games as well, but they are often less suitable for joint engagements because of the demands they make in terms of attention and rapid cognitive and physical responses. This research may be highly relevant to mathematics education settings, especially if teacher mediation is taken into account. It suggests that there may be advantages to designing environments in which teacher-student(s) conversation can be built into the technology-based activities. This may be easier to accomplish with open-ended environments in which the teacher is involved in proposing tasks or responding to students' actions, in contrast to level-driven and highly instructive environments where the parent or the teacher has little role to play.

**Digital technologies focusing on number sense**

In this section we report on studies involving digital technologies designed to support the teaching and learning of number sense [4], a fundamental aspect of early mathematical learning. As mentioned above, a wide range of digital technologies have been developed and studied, including desktop computer software, internet-based applets, touchscreen applications and, of course, calculators [5]. Some studies report on the "effectiveness", while others describe design features, or particular aspects of the "usability", that are hypothesised to support student learning. We report first on studies that involved VMs, focusing specifically on the constraints imposed and offered that contrast with the associated physical manipulative. We then consider a variety of digital tools used for different aspects of number sense while not being virtual instantiations of physical manipulatives. Finally, we present three new touchscreen applications and discuss their unique potential with respect to the development of children's number sense.

*From physical to virtual manipulatives*

Children whose learning occurs in rich environments that include (virtual) manipulatives tend to learn better and reach higher levels of academic achievement (see for example, Steen, Brookes & Lyon, 2006). However, it is not the simple presence of the (physical or virtual) manipulatives that makes the difference, but how these manipulatives are designed and used (for example, Goodwin & Highfield, 2013). Despite the abundant availability of VMs for the early years, little research has been carried out on their effectiveness and use in the classroom. VMs are, for the most part, manipulable digital technologies, both in Sedig and Sumner's sense as well as in Goodwin and Highfield's—this is not surprising given their connection to physical manipulatives. For the most part, the benefits of virtual manipulatives are seen as augmenting those of physical ones by providing more precision, more feedback that is mathematically relevant and by demanding more mathematical forms of expression (through numbers, symbols or actions).

Some recent research focuses on studying the way teachers can or could orchestrate the use of VMs [6]. As a first example, we consider the *e-pascaline*, a virtual version of the mathematical machine known as the *pascaline* (Maschietto & Soury-Lavergne, 2013; Mackrell, Maschietto & Soury-Lavergne, 2013). In the studies reported by these authors, this VM is used after students have interacted extensively with the physical *pascaline*, a mathematical machine used in many Italian primary school classrooms to foster the learning of place value (Bartolini Bussi, 2011), thereby leading to the notion of a "duo of artefacts" in which *both* physical and virtual manipulatives feature in a mathematics lessons. The *e-pascaline* has the main constitutive elements (and even colours) as the physical manipulative; however, its implementation involved additional design decisions that lead to new affordances (see Figure 1). First, all tasks are made explicit within the VM, whereas the physical *pascaline* does not include any instructions or directions. Second, the *e-pascaline*'s buttons can be hidden or shown, which affords bringing mathematical ideas and action to conscious awareness. For example, to count by 1s the tens button and the 100s button may be hidden or greyed out. Third, the wheels are turned by clicking arrows that indicate clockwise or counterclockwise rotations, thereby affording a discrete, indirect interaction that differs from the direct, continuous one of the physical pascaline. The authors explain the rationale for this decision (as opposed to, say, "click and drag to the left or right") in terms of helping children attend to the number of moves of a wheel—thus affording Sarama and Clements' sense of bringing mathematical ideas and action to conscious awareness. The authors thus make a theoretical argument that using the *e-pascaline* will help children develop a mathematical awareness that may only be left implicit with the physical counterpart. More research on how children move from these machines to paper-and-pencil forms of expressing place value may be needed in order to better understand how the *pascaline* and the *e-pascaline* function as models or analogies in children's thinking.

<Figure 1 here>

A second example of a VM is base-ten blocks, for which many applets have been designed. We highlight two specific affordances of these VMs that distinguish them from their associated physical manipulatives. The first involves the automatic transformation of a 10-block into ten individual units blocks when moved from a 10s column into a 1s column (see *Base Blocks* from the *National Library of Virtual Manipulatives*' [7] (NLVM) collection). This conversing interaction enables learners to see how the column location affects the meaning of the block while also affording Sarama and Clements' mental actions on objects. Although no empirical evidence for the effect of such a design choice has been reported, its presence in other VMs, such as Kortenkamp's *Place Value* [8] (an iPad and iPhone app), indicates some consensus about its desirability. Second, most of these VMs also show and update the numerical value of the tokens, blocks or chips placed in the different areas, thereby providing symbolic feedback and reducing the need for learners to count and calculate. This latter one seems to be an important affordance not specifically identified by Samara and Clements (2009) but potentially significant in shifting the attention of both the learner and the teacher away from computation.

Another example we include in this section is *Numberbonds* [9] (also developed for the iPad [10]), which aims to strengthen continuous and relational representations of numbers. Introducing numbers through this kind of representation emphasizes *relations* between them (a form of ordinality) as opposed to their *absolute denotation* of sets of objects (cardinality). Such an approach has been adopted, for example, in a mathematics curriculum developed by Gattegno (1970); moreover, some neuroscientific studies suggest that a more explicit and early emphasis on ordinality may be the key to learning early number (Lyons & Beilock, 2011). *Numberbonds* has a tetris-like set up in which a Cuisinaire-like rod falls in an area with a set length and the player has to quickly choose a rod from the ones displayed on the right that together with the one that has fallen completes the set base length (see Figure 2). At each level the game offers different sets of rods—until a total of 10 are displayed—to choose from to complete the reference length. Research from the fields of neuroscience and cognitive psychology conducted with students showing weak number sense highlights some advantages of this virtual adaptation of the manipulatives (Butterworth, 2011; Butterworth & Laurillard, 2010): timing of the falling allows for some "training" to occur; the changeability of the rod colour; and, the replacement of the rod's colour by a numerical value, which affords the symbolizing of mathematical concepts as well as the linking of the concrete and the symbolic through feedback. Moreover, the learning environment provides guided feedback, which enables learners to adjust their actions in relation to the goal, rather than rely on help from a teacher. Unlike the previous two examples, which replicate most of the manipulative possibilities of the associated physical manipulatives, Numberbonds addresses only one small component of the activities enabled by Cuisinaire rods; it also has a more *instructive* design, though the feedback from the environment may provide some of the mediation that features in joint engagement with media.

<Figure 2 here>

As a last example, we discuss some possible instantiations of the number line in order to introduce a more general issue of feedback. Consider a number line such as the one in *Motion Math-Fractions* [11] for the iPad, an *instructive* digital manipulative that focuses on fraction estimation. Within this environment a number line appears on the ground together with a ball that can bounce (completely elastically) and that can be controlled by the gravity accelerator of the iPad. The ball contains a fraction that has to be placed correctly on the line. Positive feedback is given for a correct placement while a hint is offered for an incorrect placement. Although this type of feedback loop may improve fraction estimation skills (Riconscente, 2013), it does not provide opportunities for learners to manipulate or express new mathematical meanings. Such possibility can be provided only through an educator's mediation (see, for example, Bartolini, Baccaglini-Frank & Ramploud, 2014). Indeed, in her study of children 5-8 years old, Goodwin (2009) found that using manipulable digital technologies to learn about fractions resulted in the production of the "most developed and advanced representations", in comparison with those using an instructive digital technology (cited in Goodwin & Highfield, 2013, p. 208). Further, these researchers note that children working with the instructive digital technology were more focused on receiving positive feedback than on discussing or reflecting on the embedded mathematical concept—as would be predicted by the joint engagement with media literature.

*Other digital tools for learning number*

We now report on studies involving digital tools that focus on different aspects of number sense while not being what we have referred to as VMs. We begin with *Building Blocks*, which is a preschool mathematics program that includes a set of different digital tools each associated with specific concepts in the curriculum and designed to meet benchmarks in a hypothesised learning trajectories (see Clements & Sarama, 2004). For example, the early levels of one of the on-computer activity sets, which are more *instructive* than *manipulative*, use a cookie-baking scenario to teach matching collections and various counting skills; at higher levels, children add. For example, a character may hide two, then one more chip under a napkin, and the child is asked to make the second cookie have the same number of chips—thus aiming at the *Non-Verbal Addition* level (Clements & Sarama, 2004, p. 183). This task could also proposed without the support of technology, but that technology provides enables the teacher (or student) to control the length of exposure of the items hidden and added; further, the system can repeat the task a large number of times, adjusting its difficulty to the input given each time by the student, and recording the performance (speed and accuracy) each time. Various other software applications have been developed to propose this kind of task, including the NCTM's "How many under the shell?" [12]. Although Clements and Sarama do not report on the gains or effects particularly related to this computer-based activity, they have shown that the curriculum as a whole increases children's development of number sense in comparison with other curricula.

In addition to classroom-oriented software programs, several remedial interventions have also been designed. For young children with difficulties with numbers [13], studies of *adaptive* and *instructive* software games such as *The Number Race* (Wilson, Revkin, Cohen, Cohen, & Dehaene, 2006, Wilson & Dehaene, 2007), *Grapho-game Maths* (Räsänen Salminen, Wilson, Aunio, & Dehaene, 2009), *NumberBonds* [14] and *Dots2Track* [15] (Butterworth & Laurillard, 2010; Butterworth, 2011) have been undertaken. *Dots2Track*, developed within the *Digital Interventions for Dyscalculia and Low Numeracy* [16] project (Butterworth & Laurillard, 2010), is based on a very simple kind of interaction in which a set of dots (or other objects) is shown on the screen and children type the corresponding number in Arabic digits, or in letters. The dots given are in fixed arrangements (for example 3 dots are in a top-left to lower-right diagonal in the stimulus area, Figure 3).

<Figure 3 here>

The software registers the response time and immediately provides corrective feedback. Preliminary results are reported as being encouraging and in the near future, the authors expect neuroimaging to be part of summative evaluation of intervention for dyscalculia. *Number Race* is another adaptive software designed to target the inherited approximate numerosity system in the IPS (Feigenson, Dehaene & Spelke, 2004) that may support early arithmetic, by presenting comparison tasks of small magnitudes presented in analogical and symbolic forms. For example, two sets of dots can appear simultaneously (and stay on the screen until the child answers, or disappear after a short time) and the child has to select the set with greater numerosity (Fig. 4a). Wilson and colleagues (Wilson, Dehaene, Pinel, Revkin, Cohen & Cohen, 2006) show how the software is designed to create a multidimensional model of the evolution of a learner's "knowledge space" along three variables and time (Fig 4b). After five weeks, experimental subjects showed greater improvement in their number sense with respect to the control group (Wilson & Dehaene, 2007).

<Figure 4 here>

*Potentials of multitouch technology*

The previous examples in this section have been designed for the mouse and keyboard inputs available with desktop and laptop computers. The main interaction is through clicking (rather than typing) in discrete objects, with one child working with the software at a time. In contrast, the following three examples are prototypical in exploiting various potentials of multitouch technology with respect to number sense learning. With multitouch technology, the interaction becomes more immediate, as the fingers contact the screen directly, either through tapping or a wide variety of gestures. Further, the screen can be touched by multiple users simultaneously at the same time, which invites different types of activity structures than the computer or laptop. We have devoted a section to these new technologies because of the close link between their main interaction mode (through the fingers) and the emerging neuroscientific literature pointing to the importance of fingers in the development of number sense (e.g., Butterworth, 1999a, 1999b; 2005; Penner-Wilger, Fast, LeFevre, Smith-Chant, Skwarchuk, Kamawar, & Bisanz , 2007; Gracia-Baffaluy & Noël, 2008; Andres, Seron, & Olivier 2007; Sato, Cattaneo, Rizzolatti, & Gallese, 2007; Thompson, Abbott, Wheaton, Syngeniotis, & Puce, 2004). Basic component abilities that can be powerfully mediated through multi-touch technology are: 1) subitising; 2) one-to-one correspondence between numerosities in analogical form and fingers placed on screen/raised simultaneously/counting with fingers, and in general finger gnosia; 3) fine motor abilities; and, 4) the part-whole concept. While the examples we discuss here may quickly be replaced with newer versions, they allow us to identify particular design decisions that will be relevant to new applications as well.

Our first example is an *instructive, conversing* digital technology called *Fingu* [17] (Barendregt Lindstrom, Rietz-Leppanen, Holgersson, & Ottosson, 2012), an iPad application in which the stimuli are given as fixed arrangements of floating objects (see Figure 5a). Users must place the corresponding number of fingers on the screen simultaneously and each touch produces a fingerprint—thus providing visual feedback and rapid closure to the *conversation* (see Figure 5b), which ends with an additional yes/no feedback (auditory signal as well as visual cues: happy animations for correct responses and sad animations for incorrect ones). The game is timed and at each of the seven levels the number of objects that appear increases while the time to respond decreases. Very little opportunity for joint engagement is offered.

<Figure 5 here>

The necessity for simultaneous rather than sequential input further encourages subitising. The choice of floating groups of objects (where the disposition of each group remains invariant) differs from the fixed arrays typically offered in other environments, provides a different type of stimulus for the solver's object tracking system (Fayol & Seron, 2005; van Herwegen, Ansari, Xu, & Karmiloff-Smith, 2008; Cantlon Safford, & Brannon, 2010), a system upon which the ability to subitise supposedly relies (Piazza, 2010). In an exploratory micro-longitudinal pilot study with 5- and 6- year old students, researchers found that children playing the game for a three-week period, with guidance from the teacher, would score significantly lower in higher levels of the game, but overall would increase their percentage of right answers (Barendregt *et al.*, 2012). The researchers also identified some indicators (counting all, counting from smallest, counting from largest, counting fingers to 5, counting fingers over 5, manipulating fingers with other hand, problem pressing down) that could contribute to understanding different learning trajectories in *Fingu* that may lead to improvement of children's mathematical abilities. This is claimed as one of the goals of an on-going study (Barendregt *et al.*, 2012, p. 4).

Ladel and Kortenkamp (2011) report on an open-ended, *manipulable, constructive* multi-touch environment developed to foster children's development of the *part-whole concept* [18] (as proposed, for example, by Resnick, Bill, Lesgold & Leer, 1991) and of *finger symbols sets*. These are thought to foster flexible calculation strategies, such as composing and decomposing numbers with respect to 5 and 10 (Brissiaud, 1992). The environment consists of a table connected to a computer that recognizes multi-touch inputs. Children can create and move tokens simultaneously (through multi-touch but also multi-user). Tokens can also be "fused together" (for example 3 and 5 to make 8) and the environment will give symbolic feedback (3+5=8) on the action. It is also possible to give tasks in symbolic form (e.g. "3 + 4 = _" or "3 + _ = 7"), which the children are asked to express with tokens. Ladel and Kortenkamp (2012) report on a different task in this environment designed for children (age 5-7), in which a certain number of virtual tokens must be placed on the table "all at once", thus focussing on the shift from sequential, ordinal counting to holistic cardinal count [19]. The researchers describe the different ways in which the children attempted to solve the task, stressing how the environment enabled the children to exteriorise their thinking about number, and noting the important role of the teacher in interpreting and responding to the children's actions.

*TouchCounts* (Sinclair & Jackiw, 2011) is made up of two constructive microworlds involving *manipulable*, *constructive* interactions: ("1, 2, 3,…") and ("1+2+3+…"). Audio feedback can be given in English, Italian or French. In the "1, 2, 3,…" microworld, the mode "with gravity" presents a "shelf" on the screen (Figure 6a). As a finger is placed on the screen, a coloured circle containing an Arabic numeral appears on the screen and the number is also spoken orally. When the finger is lifted from the screen the numbered object falls, unless it is dragged so that it "sits" on the shelf (see Figure 6b). When more fingers are placed on the screen, the counting continues on from the last number reached with the previous touch. If the child interacts through successive touches using a single finger, or different fingers placed on the screen sequentially, *TouchCounts* will end up "counting" for the child, so it may strengthen her recollection of the number words (in sequence) and possibly fine motor abilities. If multiple fingers touch the screen, the same number of numbered objects appears but only the last number is given orally—this enables, for example, counting by twos. In the "no gravity" mode, the numbered objects do not fall so that placing five fingers on the screen (sequentially or simultaneously) will produce something similar to Figure 6c.

<Figure 6 here>

This type of interaction involves various number sense component abilities such as subitizing and fine motor skills (simultaneous touch to generate the proper enumeration and dragging of a selected finger to place the circle on the line), but the environment may also help lay foundations for the counting principles and for the transition from ordinal to cardinal counting. This is because it can foster the development of awareness of the one-to-one correspondence between fingers and numerosities, or between numbers and successive touch-actions on the screen (one-to-one correspondence principle); it may foster memorization of the sequence of number-words to recite when counting (stable order principle); in the modality without gravity, the last word heard through the audio feedback corresponds to the total number of circles on the screen (cardinality principle); and finally in the gravity mode the possibility of marking certain numbers by dragging them on the line may favour a process of reification (Sfard, 2008) of the number, necessary for operating on numbers.

In the "*1+2+3…*" microworld, touching the screen with several fingers simultaneously generates sets of circles, of a same random colour, enclosed in a numbered circle indicating its magnitude (Figures 7a, b). Through the "pinch" gesture (Figure 7c) it is possible to act on the cardinal numbers, in this case, adding them together. This gesture constitutes a fundamental metaphor of addition, that of "grouping together" (Lakoff & Núñez, 2000). The gesture can be performed prior to any formal instruction about addition, of course, but may help children develop a metaphoric meaning for addition, as well as a sense of the symmetry of this operation (and thus the commutative property). The result of pinching two groups together is a larger group in which the colours of the addend circles have been preserved in order to leave a trace of the action (Figure 7d).

<Figure 7 here>

Pilot studies indicate that with the use of certain tasks, children 5-6 years old can learn to shift from thinking of number in terms of the process of counting to thinking of them as reified objects (Sinclair & Heyd-Metzuyanim, 2014). Further, with even younger children (3-4) years old, pilot research has shown that the task of placing fingers all-at-once on the screen can help develop their "gestural subitizing" (Sinclair, 2013). These studies show how a rich *manipulative* and *constructive* environment can be also be *instructive*, thus fostering the development of both procedures and associated concepts.

While research is in its early phases with respect to new touchscreen and multitouch environments, we propose the following summary of the basic features of multi-touch technology that can be used to foster the development of number-based concepts and abilities.

<Table 1 here>

As has been the case with other significantly new digital technologies, we also expect to see some changes in the way number sense concepts themselves may develop—along with the order and pace by which these concepts are learned—as children interact not only through alphanumeric means but through touch, sound and image as well. While the affordances seem to be clearly geared to supporting young children's learning, more research is needed into how these touchscreen devices might affect children's interactions with physical tasks involving number sense, including the still pervasive pencil-and-paper media of the mathematics classroom.

**Digital technologies focused on geometry**

The goals of geometry learning at the k-2 level are to develop a better understanding of objects in relation to their shape and position, and to attend to the geometric properties (parallelism, congruence, symmetry) and behaviours (invariance, sameness) of these objects. Since similar goals pertain also the higher elementary and middle school geometry, some of the research already conducted at these later grades using digital technologies is relevant to this chapter [18]. Dynamic geometry environments (DGEs) have been the most widely researched geometry technologies and researchers have shown how they help foster conjecturing, enable learners to interact with a larger number of examples, and help learners attend to invariances (see Baccaglini-Frank & Mariotti, 2010; Laborde, Kynigos, Hollebrands and Sträßer, 2006; Mariotti, 2006; Sinclair & Robutti, 2013) in a wide range of geometry-related topics. Many of these findings relate to the dragging tool available in DGEs, which enable direct, continuous manipulation (and this form interaction is now also available in some web-base applets and VMs aimed at young children). For example, Battista (2007) has shown that the dragging tool, used to transform constructed quadrilaterals, enabled grades four and five children attend to the invariant properties of the different quadrilaterals and even identify certain quadrilaterals (rectangles) as subclasses of others (parallelograms). At the high school level, Hollebrands (2003) has shown that the use of dragging can help students develop a functional understanding of transformations (reflections and rotations) since dragging one object continuously on the screen changes the position of another associated object.

*Example from the Building Blocks curriculum*

In addition to the *Building Block* VMs mentioned above, Clements & Sarama (2002) describe a geometric *manipulable* and *constructive* digital technology called "Piece Puzzler". It was intentionally designed to contain screen versions of pattern blocks and tangram shapes, which children can manipulate to create or duplicate larger composite shapes. The authors report research results on the effect of the curriculum as a whole (Clements & Sarama, 2007), but not on the specific use of the virtual manipulative. However, the principles used for its design are of interest for several reasons. One such principle stipulates that the virtual manipulative, along with the specific tasks, be designed in terms of a hypothetical learning trajectory. In this case, the authors propose seven levels along this trajectory, each of which is accompanied by a specific task aimed at achieving the particular level, with the final level aiming at having children be able to iteratively compose composite shapes to tile the plane. For example, the goal of level four (Shape Composer) is for children to choose and manipulate, through turning and flipping, given shapes to completely fill a region. The given regions are multiply cornered so that children have to attend to angles in the shapes as well as in the region [20]. In level 6, the children can use shapes to create a composite objects in the shape of a toy (like a rocket ship), which they can then duplicate by pressing the "do it again" tools, thus creating iteratively composite units. In the final level seven, children create superordinate units of tetrominoes and use them to tile the plane.

In their early reporting on the use of the Building Blocks program, Sarama and Clements (2002) report that the "use of the tools encourages children to become explicitly aware of the actions they perform on the shapes" (p. 103) since, unlike physical pattern blocks and tangram shapes, children cannot just pick up and move the pieces with their hands. Further, the children are able to create designs that are more precisely assembled than if they were working with physical objects since, as Moyer, Niezgoda and Stanley (2005) have also pointed out in their study of the NLVM Pattern Block, the shapes can be "snapped" into position, and stay fixed. This description of the software hides many design choices that are central in determining how children use it and what they learn as they use it. Clements and Sarama (2002) acknowledge that their pilot testing raised questions about whether unexpected outcomes should lead to changes in the software design or changes in the learning trajectory. This should not be surprising if one acknowledges how tools mediate learning, which is in keeping with Clements and Sarama's general Vygotskian approach.

In terms of design, the Clements and Sarama briefly discuss decisions made around how the turning would be handled by the software. They tested four possible choices (tool, button, direct manipulation with continuous motion and direct manipulation with discrete units). Their ultimate choice of a tool interface, which is discrete and direct) was motivated by the fact that three and four year old children found this interface easier to learn and to use. However, it requires the designer to choose a default turn angle, which means that children's turn actions are actually 'turn-by-30-degrees' actions (for the pattern blocks, and fifteen degrees for the tangrams), thereby highly constraining the example space of turn (not to mention the fact that the centre of rotation is always in the middle of the shape). Similarly, the flip tool makes a choice that the line of reflection will be horizontal and immediately under the shape being flipped. Given that the intended learning trajectory, these constraints on reflecting and rotating may not be too problematic. However, they are the kinds of constraints that may no longer be needed with touchscreen interfaces, where children can act on objects with their fingers instead of through the intermediary of the mouse. Indeed, in the next section, we describe a learning environment in which the turn interface chosen is that of direct manipulation with continuous motion.

An example that focuses on very similar mathematical ideas, but that differs both in the design of the interface and the accompanying task, comes from a project involving the use of *Cabri Elem*, which is a multi-purpose dynamic geometry environment that can be used to design microworlds suitable for primary school learners. The "Tiling" microworld, which is both *manipulable* and *constructive*, involves the composition of shapes into tiling patterns. Children can manipulate eight different shapes directly by sliding, turning or flipping. Laborde and Marcheteau (2009) report on a study conducted with grade three children, who were asked to work in pairs to create tilings involving at least two different shapes and then to describe their pattern to another pair who would try to re-create the pattern. While all fourteen pairs successfully created tilings, only three were able to describe these tilings in terms of a repeating unit. Given the success of the Clements and Sarama's *Piece Puzzler*, as well as the NLVM *Pattern Block* used by both Moyer *et al.* (2005) and Highfield and Mulligan (2007) with younger children, the "cloning" affordance through which a user can duplicate an existing set of shapes, seems instrumental in promoting thinking in terms of repeating units.

*Example of whole classroom dynamic geometry*

In the third example, we offer a cluster of examples that focus on different aspects of geometry at the k-2 level, and that involve a plenary mode of interaction in the classroom with a DGE being projected on a wall or an interactive whiteboard. The first example concerns Sinclair and Moss's (2012) study involving kindergarten (4-5 years old) children and triangle identification in which *Sketchpad* is used as a *conversing* and *manipulable*, as well as *constructive* technology. In this study, which used Sfard's communicational approach in which learning is conceptualised as a change in discourse, children moved relatively quickly from a first van Hiele level (in communicational terms, a discourse about the physical reality around us where shapes are identified as the same if they match) to a second one (in communicational terms, a discourse that treats level one things as objects and where shapes are identified as the same if they can be transformed one into the other), with some even moving to a third van Hiele level (in communicational terms, a reified one, where shapes are identified as the same by comparing verbal descriptions of the shape). In particular, the children initially used the word *triangle* much like a proper name corresponding to an equilateral triangle and identified shapes as triangles if they looked like an equilateral triangle. When these children were shown a triangle with its vertex pointing down, which was constructed using the segment tool in *Sketchpad*, they either said that it was an "upside down triangle" or that they could see it as a triangle if they turned themselves upside down. After the teacher dragged one of the vertices of the triangle on the screen, all but one child began to speak of triangle as a family name that describes a larger set of triangles than just equilateral ones. They identified these non-equilateral triangles using a routine of transformation in which a shape was a triangle because you could "stretch it out." A few children even began describing the triangles on the screen in definitional ways, stating that "[e]very triangle could be, um, a different shape but it just has three corners" (p. 36). The children were given the opportunity to create their own triangles using the teacher's computer and made a variety of triangle shapes, including long and pointy ones and upside-down ones. When the teacher shifted to squares instead of triangles, some of the children immediately used a discourse of transformation when talking about a square that was sitting on its vertex, while others insisted that this shape was a diamond. This suggests that the children had not all succeeded in shifting discourse, but this is not surprising given the short intervention (30 minutes).

A more extended intervention described in Sinclair and Kaur (2011) involves introducing young children (kindergarten and grade one) to the concept of angle, this time using *Sketchpad* in a *manipulative/manipulable* manner. While angle is typically not formally included in the curriculum at this age, the researchers deemed it both possible and desirable to introduce it in a visual and dynamic way rather than a measurement-based one (which involves learning about degrees and using relatively large numbers like 180 degrees). The goal of the study was to determine whether focusing on the metaphor angle-as-turn might enable young children to develop understandings about angle and help address common errors that students make, as identified in the literature, such as assuming that an angle with longer arms is bigger than an angle with shorter arms. The initial task in the intervention focuses on developing benchmark angles. This is done by offering children clockwise and counter-clockwise turning options of ¼, ½, ¾ and 1—when one of these buttons are pressed, the arrow turns dynamically and leaves a trace of the swept out angle (see Figure 8). The teacher prompts the children to arrive at a given destination using a smaller or bigger angle so that the word "angle" comes to be associated with the turning motion. Then, children are invited to use the sketch shown in Figure 8 with the task of getting the car to the gas station. By dragging the angle "dial", children determine the amount by which the car will turn. This requires children to coordinate the turning of the dial with the turning of the car, which is challenging when the car is not facing up, like the vertical arm of the angle dial. The disassociation of the car form the angle dial is thus crucial to the design of the sketch.

<Figure 8 here>

The researchers found that the children had no difficulty using benchmark angles to describe the position of objects in the screen. Further, although the children initially struggled with the angle dial, they were eventually able to use it to get the car to the gas station as well as to solve more complicated, multi-step trips. Follow-up computer-based as well as paper-and-pencil based tasks showed that the children, when asked to compare angles (with distracting features such as differently-sized arms and different orientations) focused on the amount of turn rather than the size of the arms. Further, when working on tasks involving triangles, the children were as likely to focus on the number and size of the angles than on the sides. The dynamic nature of the software, along with the immediate, non-evaluative feedback, enables children to see the turning action and experiment with the effect of different sizes of turns, while the teacher's questions and interactions help the children associate the word 'angle' to this turning action. This research offers initial evidence to support the claim that angles could be effectively introduced earlier in the geometry curriculum, which would have implications for the way other topics are taught.

The final example involves reflectional symmetry. Ng and Sinclair (to appear) describe a three-lesson intervention involving both computer-based and paper-and-pencil based activities aimed at helping children attend to the geometric properties of reflectional symmetry. In the first lesson, the children explore the *manipulative/manipulable* "symmetry machine" through direct, discrete interaction by dragging one of the squares on the screen and observing the resulting motion of the associated square (see Figure 9a, b). The discrete motion of the square is meant to draw students' attention to how the movement of the square and its image are related. After initial exploration, the children are shown various designs and asked whether they could be created using the symmetry machine. If so, the students recreate the design, with the environment providing self-checking feedback, and if not, they are invited to explain why. The children can also interact through direct, continuous dragging with the line of symmetry in order to reproduce designs that have horizontal or oblique, as well as vertical symmetry (see Figure 9c). In the second lesson, children are shown a broken symmetry machine containing squares only on one side of the line of symmetry and they are asked to fix the other side to make it symmetric. On the third lesson, the children use the continuous symmetry machine, which involves direct, continuous dragging of a traced point, as does a symmetric point. The children are asked to create certain shapes using this continuous symmetry machine (butterfly, a square, a house, etc.).

<Figure 9 here>

The researchers found that the children developed new and embodied ways of thinking about symmetry and its properties. Further, the children moved from a static conception of symmetry to a functional and dynamic one, focusing on the symmetric relation between a shape and its image rather than on the static property of being symmetric. This shift was occasioned by the processes of semiotic mediation in which the dragging and tracing tool, as well as the language and gestures of the teacher, became signs that enabled communication about central features of reflectional symmetry including: the way in which one side of a symmetric design is the same as the other; the way in which one component of a symmetric design is the same distance away from the line of symmetry as its corresponding image; the way in which a pre-imagine component and its image have to be on the same line relative to the line of symmetry; and, the way in which a pre-image and an image gives rise to parity.

**Learning mathematics through programming**

Programming is an important topic that some countries, such as Italy and the UK [21], even include explicitly in the curriculum indications for preschool or primary school [22]. However, after the research on Logo, little research has been published on the teaching and learning of mathematics through computer programming. Currently, a variety of programmable toys (such as Bee-bot, Probot and Lego NXT) are being used in classrooms around the world, as may be seen by the number of activities available (especially for bee-bot) on websites from different countries, and the use of such toys is proposed in textbooks for college courses for Italian pre-service teachers (see Baccaglini-Frank, Ramploud & Bartolini Bussi, 2012) or in special documents for Australian in-service teachers [23]. However as reported by Highfield and Mulligan, "The consequence of young children's immersion in these technologies has not been adequately investigated and the potential advantages and disadvantages for mathematics education need to be examined" (2008, p. 1). In a study aimed at comparing mathematical processes that preschool children might encounter whilst 'playing' with different forms of technology, including the bee-bot, the researchers report that it remains "unclear how mathematical processes are explored, understood or assimilated by young learners" (2008, p. 6), and teacher assistance and guidance is heavily needed. A similar finding is reported by Pekarova (2008), who notes that "some children clearly demonstrated deep comprehension of principles of Bee-Bot's control" (p. 120). These findings resonate with what is reported in theses of graduate students at the Department of Education and Human Sciences at the University of Modena and Reggio Emilia (Italy) who have carried out a number of preschool and early elementary school interventions using the bee-bot (e.g. Bartolini Bussi & Baccaglini-Frank, 2015). Other studies propose the use of robotic toys to foster problem solving, mapping, and measurement activities (Highfield, 2010; Highfield, 2009; Highfield & Mulligan, 2009).

One of the few published studies involving the use of a programming environment at the k-2 level is that of Yelland (2002), who developed a computer microworld based on the use of *Geo Logo*, which encouraged young children (with a mean age 7 years and 4 months) to explore concepts of length measurement. The microworld enabled the children to move different turtles, each of which took a different length of step. The children worked in pairs and had to decide how to coordinate the length of each step with the number of steps in order to reach a destination. Yelland found that the computer-based environment was more conducive to interactions between the pairs of children "which forced them to use number and compare numbers in new and dynamic ways" (p. 86). She further opined that the environment "facilitated playing with units of measurement in ways that were not possible without the technology" (p. 91). These findings are consistent with research conducted at the higher-grade levels (see, for example, Clements & Battista, 1989). It also suggests that children at this age are capable of more sophisticated reasoning with measurement—in terms of co-varying the length of the step and the number of steps—than current curricula assume (also see Goodwin & Highfield, 2013; English & Mulligan, 2013).

**Discussion**

The rate of proliferation of new digital media far outpaces the amount of research that can be conducted. This is in part due to the ever-changing developments in hardware and software, and also to the gap between the kinds of environments valued by researchers and those produced by software design companies. Indeed, several of the examples we have analysed have been designed by the researchers themselves, who generally favour more open-ended environments that support both conceptual and procedural aspects of concept development. However, funding for such software development can be scarce, and support for ongoing maintenance even more so. Software is not like textbooks and needs continuous debugging, upgrading and adaptation to new operating systems and hardware affordances. Although it is likely that many of the environments discussed here will be superseded in the years to come, our analyses should be understood from a broader perspective, and a reader should appreciate how they will remain the same regardless of the "game" or "interface". To that end, we have highlighted some of the different design decisions that are relevant to mathematics learning at the k-2 grades.

Multitouch devices seem to have particularly strong potential in children' development of number sense, in large part because of the important role that fingers play in this development and also because of the direct mediation enabled by this technology. Other computer-based programs, including VMs, have also been shown to be effective, especially when the feedback—either from the computer or a teacher—can help children reflect on the mathematically relevant aspects of their actions. We expect that many of these VMs will soon be available on multitouch platforms as well, but the specific ways in which they are implemented—in terms of the kinds of interactions that are made possible—will be important in determining their effectiveness. In terms of geometry, a large portion of the studies have focused on microworlds that have been built within dynamic geometry environments. The ease with which young children can see and explore a wide variety of continuously changing shapes, as well as relationships between shapes, seems to make them highly suitable for the k-2 level, particularly with appropriate teacher mediation. We expect to see more standalone microworlds developed for multitouch platforms, which will lead to interesting research questions about the different activity structures that would enable several children at once to interact with shapes, perhaps working together to compose a set of polygons.

Clearly, much more research on the effects of digital technologies on k-2 teaching and learning is needed. Although there exists an abundance of digital tools for researchers and teachers to choose from, especially in terms of internet applets and touchscreens apps, we know very little about how particular design choices might affect children's learning—as well as children's use of physical materials in the classroom. We would like to draw attention to specific themes that we think are and will become significant in this area of research, including: choosing between discrete and continuous flow of interaction; accounting for the affective dimension of learning; acknowledging the impact of different types of feedback; and, attending to constraints.

While earlier software programs often chose more discrete modes of interaction (as in Clements & Sarama's discrete angle interface and the NVLM's discrete number line applets), there may be a move toward more continuous models both because of the new forms of interactions afforded by dynamic and touchscreen technologies and because of changes in our assumptions about what mathematical objects children should be exposed to. In terms of the latter, if prior technologies privileged whole numbers (Cuisinaire rods, base ten blocks, and even discrete numberlines), new technologies enable and sometimes require children to work with real numbers (Crespo & Sinclair, 2006), which presents interesting challenges for VMs as well as, more generally, for orchestrations of digital and physical classroom resources. We are reminded of Papert's adage promoting a "Mathland" learning environment that, much like living in a foreign country, does not limit forms of participation and interaction to the novice's basic vocabulary.

Many of the research studies we read, as well as reports in professional journals, highlight the higher level of engagement that children seem to experience when working with digital technologies. Although we assume this is not universally true, it is also part of the anecdotal evidence surrounding digital technology use at the high grades. Cognizant of the negative effects that too much "screen time" can have on young children (see Public Health England, 2013), more research on the nature and function of that engagement is warranted. It is not always clear from research reports whether the enjoyment stems from a simple change in environment or from the tighter feedback loop that the digital tool provides or from the kind of curiosity, pride and intellectual engagement that is typically associated with self-directed and deep mathematical learning. This would require theoretical and methodological tools that enable analyses that do not dichotomise the cognitive, affect and aesthetic dimensions of learning.

One of the frequently-mentioned strengths of digital technologies is their ability to provide instantaneous, customized feedback, which is very challenging for a classroom teacher to do. However, as we have signaled throughout the chapter, there are a wide variety of forms of feedback, each potentially functioning very differently both cognitively and affectively for children. More work is needed to understand how young children process and use these different forms of feedback [24]. For example, given the extensive research showing that evaluative, corrective teacher feedback (of the yes/no form) in the classroom has negative effects on student learning, developing and using digital technologies that only have this kind of feedback must have a strong counter-rationale. At the same time, for teachers and research opting for more guided forms of feedback (responding through hints, for example), research that studies the way that children use—and what they learn from—this feedback would be very helpful. In addition, the exteriorization of thinking offered by digital technologies (which Noss & Hoyles (1996; 2006)) refer to as a "window on mathematical meanings") can be extremely helpful to teachers. This implies that more research is needed on how teachers make sense of and use feedback in working with and assessing their students.

We offered two different frameworks for distinguishing various digital technologies, both of which provided comparative power, but neither of which helps designers, researchers and teachers attend to the way particular constraints can affect the development of mathematical meanings. Some constraints are determined by the technology (as in Clements and Sarama's choice of a discrete interaction) while others are the result of didactic goals. Constraints can be used purposefully to restrict and focus a learner's interactions, as discussed by Ladel and Kortenkamp (2011; 2012) and by Manches, O'Malley and Benford (2010). For example, a constraint in *TouchCounts* is that numbered objects cannot be moved once they are placed on the screen; this constraint aims to restrict the children's externalizing actions in order to the internalization of the one-to-one-to-one correspondence between touch, objects and number. In order to fully understand the potential of digital technologies, and their impact on student learning, the nature and consequence of these constraints needs to be taken into account—this will require better communication between designers, researchers and teachers in future research. Further, methodological choices for research need to reflect the way in which the learning environment as a whole, including the teacher and the tasks used in concert with the digital technology, functions to affect mathematical learning. To this end, we find persuasive Stylianides and Stylianides' (2013) argument in favour of teaching experiment (or classroom-based interventions) methodologies, which can increase the likelihood that the results of research are applicable while also shedding light on how and why certain situations work.

**Endnotes**

1 For example, dynamic geometry environments (such as *The Geometer's Sketchpad* and *Cabri-Elem*) can be used to create microworlds that may function as virtual manipulatives or used on interactive whiteboards (with Sketchpad also available for iPads).

2 Sedig & Sumner (2006) distinguish between "real" and "perceived" affordances as a way of underscoring the fact that users do not always use digital technologies in the ways intended by their designers. This fact is central to the theory of instrumental genesis, which is discussed in the next section.

3 Sedig and Sumner (2006) point out that such interactions are often seen as reducing the cognitive load on learners, especially if they do not have to plan their actions in advance, which they might be more inclined to do in discrete, conversing interactions. The choice of a manipulating interaction over a conversing one will depend on the goal of the tool/task/concept.

4 We acknowledge that there are some major differences between the ways in which number sense is defined in the mathematical cognition literature and its definition in the literature in mathematics education (see, for ex., Berch, 2005). Here we will refer to the literature in both fields, assuming that the more specific meaning within the former field is a necessary stepping stone towards its vaster, much more complex and multifaceted connotation proposed in the latter.

5 There is extensive literature showing the value of using calculators in the early years (Groves & Stacey, 1998) However, there seems to be much less enthusiasm for them amongst teachers and researchers, perhaps in part because they do not offer the visual forms of interaction that other digital technologies at this early age now do.

6 Research on the use of digital technologies with younger learners seems to have coincided with a new focus on the role of the teacher in technology-based classroom environments—this focus having been prompted by the realization that despite high levels of accessibility and institutional support, as well as supporting research, the use of digital technologies in the mathematics classroom remains relatively low (see Laborde, 2008).

7 See http://nlvm.usu.edu/en/nav/topic_t_1.html

8 See https://itunes.apple.com/us/app/place-value-chart/id568750442?mt=8

9 To play the game visit http://number-sense.co.uk/numberbonds/

10 See https://itunes.apple.com/us/app/number-bonds-by-thinkout/id494521339?mt=8

11 See http://www.appstore.com/motionmathhd

12 See From the NCTM Illuminations website: http://illuminations.nctm.org/Activity.aspx?id=3566

13 These children may be referred to in the literature as being affected by *dyscalculia* (Butterworth, 2005).

14 To play the game visit http://number-sense.co.uk/numberbonds/

15 To play the game visit http://number-sense.co.uk/dots2track/

16 See http://low-numeracy.ning.com/ and http://number-sense.co.uk/

17 See https://itunes.apple.com/it/app/fingu/id449815506?mt=8

18 Something analogous to the part-whole concept is present in the Chinese mathematics education tradition in its presentation of "variation problems" (Sun, 2011).

19 Indeed, they are also relevant for topics other than geometry. However, it would seem that geometric environments such as DGEs have been more amenable to stretched across the grade levels.

20 In a study focused more specifically on children's patterning, Moyer *et al.* (2005) also studied the use of a Pattern Block VM with kindergarten children and found that the patterns they created were more creative, complex and prolific when using the VM than when using concrete materials. Similarly, Highfield and Mulligan (2007) report that preschool children using a Pattern Block VM as well as a drawing tool called *Kidpix* experimented with more patterns, created more precise patterns and made more use of transformations than children who worked only with physical materials. The authors do caution that the children found the use of the mouse challenging—again, this is a hardware limitation that touchscreen technology can mitigate—and the additional affordances of *Kidpix* sometimes distracting.

21 For example, in Italy, the national curriculum stipulates that by grade three, when possible, students should be introduced to some programming languages that are simple and versatile in order to develop a taste for the planning and realization of projects and in order to understand the relation between coding language and its visual output. In the UK national curriculum, children are expected to develop competence in two or more programming languages by the age of eleven.

22 We highlight several new programming languages are specifically been developed for young children, such as *Squeak, ToonTalk* and Giorgi & Baccaglini-Frank's (2011) *Mak-Trace* app for the iPad.

23 See http://www.learningplace.com.au/deliver/content.asp?pid=38840.

24 It seems like not all children do this in the same way. See Goodwin and Highfield (2013, p. 214).